\documentclass{amsart}

\usepackage{amsfonts, amsmath, amsthm, amscd}

\newtheorem{theorem}{Theorem}
\newtheorem{lemma}[theorem]{Lemma}
\newtheorem{defn}[theorem]{Definition}
\newtheorem{prop}[theorem]{Proposition}
\newtheorem{cor}[theorem]{Corollary}

\newcommand{\D}{\displaystyle}
\newcommand{\R}{\mathbb{R}}
\newcommand{\C}{\mathbb{C}}

\newcommand{\HH}{\mathcal{H}}
\newcommand{\T}{\mathcal{T}}
\newcommand{\B}{\mathcal{B}}
\newcommand{\del}{\bar{\partial}_E}
\newcommand{\tri}{(\del, \theta, \phi)}
\newcommand{\End}{\textrm{End }E}
\newcommand{\hatb}{\hat{\B}}

\begin{document}

\title{Birational Equivalence of Higgs Moduli}
\author{Mridul Mehta}

\curraddr{Max Planck Institute for Mathematics in the Sciences,
Inselstrasse 22, D-04103 Leipzig, Germany}

\email{mehta@mis.mpg.de}

\subjclass{14D20}

\begin{abstract}
In this paper, we study triples of the form $(E, \theta, \phi)$ over a compact 
Riemann Surface, where $(E, \theta)$ is a Higgs bundle and $\phi$ is a global 
holomorphic section of the Higgs bundle. Our main result is an 
description of a birational equivalence which relates geometrically the moduli 
space of Higgs bundles of rank $r$ and degree $d$ to the moduli space of Higgs 
bundles of rank $r-1$ and degree $d$.
\end{abstract}

\maketitle

\section{Introduction}
\setcounter{theorem}{0}
\setcounter{equation}{0}

In \cite{Br1}, Bradlow introduced and studied the moduli space of pairs of the 
form $(E, \phi)$ over a compact K\"ahler manifold $X$, where $E$ is a 
holomorphic vector bundle of fixed rank $r$ and degree $d$ over $X$ and $\phi$ 
is a holomorphic section of $E$.
The construction of the moduli space of these pairs involves a 
choice of linearization, which results in the notion of parameter dependent 
stability for these objects. This parameter is a real number taking values in 
the closed interval $\left[\frac{d}{r}, \frac{d}{r-1}\right]$. Consequently, the 
construction yields not one but a family of moduli spaces $\mathcal{B}_{\tau}$ 
of `$\tau$-stable' pairs. It was shown that for all but finitely many values of 
$\tau$ in the specified interval, the spaces $\mathcal{B}_\tau$ are birational. 
Further, when $X$ is a Riemann surface, the $\tau$ stability condition 
forces the spaces $\mathcal{B}_{\frac{d}{r}}$ and $\mathcal{B}_{\frac{d}{r-1}}$ 
to be closely related to the moduli spaces $\mathcal{M}(r,d)$ and 
$\mathcal{M}(r-1,d)$ of vector bundles of rank $r$, degree $d$ and rank $r-1$, 
degree $d$, respectively. 
As a result, the birational equivalence of the moduli spaces $\mathcal{B}_\tau$ 
relates geometrically the moduli spaces $\mathcal{M}(r,d)$ and $\mathcal{M}(r-1, 
d)$ over $X$. This setup, which was used by 
Thaddeus in his famous paper (\cite{Th1}) provides a natural method to study the 
moduli space of vector bundles over a Riemann surface inductively using rank.

This paper provides a similar setup to study the moduli space of Higgs bundles. 
We consider objects of the form $(E, \theta, \phi)$ over a Riemann surface $X$, 
where $E$ is a holomorphic vector bundle of rank $r$ and degree $d$ over $X$, 
$\theta$ is a Higgs field and $\phi$ is a holomorphic global section of $E$.
Constructing the moduli space for these objects leads to a notion of 
parameter dependent stability, where the parameter $\tau$ is a real number which 
takes values in the interval $\left[\frac{d}{r}, \frac{d}{r-1}\right]$. We show 
that for 
all but finitely many values of $\tau$, these moduli spaces are birationally 
equivalent. Finally, we relate the moduli spaces of these triples at the end 
points of the interval to the moduli spaces $\mathcal{H}(r,d)$ and 
$\mathcal{H}(r-1, d)$ of $GL_r(\C)$-Higgs bundles over $X$, which completes the 
setup. The main results of this paper are the following.

\begin{theorem}
For all noncritical values of $\tau$ in $\left(\frac{d}{r}, \frac{d}{r-1} 
\right)$, the spaces $\B_\tau$ of $\tau$-stable Higgs triples are all 
birational.
\end{theorem}

Let $\HH_0(r,d) = \left\{ (\del, \theta) \in \HH(r,d)\,
\Big|\, \textrm{det}(\theta) = 0 \right\}.$ We then have:

\begin{theorem}
There exists a fibre space $\mathcal{F}$ over $\HH_0(r,d)$ which is birational 
to a projective bundle $\mathcal{E}$ over $\HH(r-1,d)$.
\end{theorem}

Finally, it should be mentioned that much of this work has been motivated by
the paper on stable pairs by Bradlow, Daskalopoulos and Wentworth \cite{BDW}.

\section{Preliminaries}
\setcounter{theorem}{0}
\setcounter{equation}{0}

We begin by fixing our notation, which we shall use for the rest of this paper. 
Let $X$ be a compact Riemann surface of genus $g$. We fix a K\"ahler form 
$\omega$ on $X$, normalized so that $\int_X \omega = \textrm{Vol}(X) = 4\pi$. 
Let $E$ be a fixed smooth ($C^{\infty}$) complex vector bundle on $X$ of rank 
$r$ and degree $d$. We will assume that $r$ and $d$ are 
coprime. We shall denote by $\mathcal{A}$ the affine space of all $\del$ operators
on $E$.

In order to be rigorous, we use a similar convention as Atiyah and Bott 
(\cite{AB}), and for 
any Hermitian bundle $V$ over $X$, by $\Omega^{p,q}(V)$ we will always mean the 
Banach space of sections of Sobolev class $L^2_{k-p-q}$, of the bundle of 
differential forms of type $(p,q)$ with values in $V$ (any $k\ge 2$ suffices).
Then $\mathcal{A}$ is the space of holomorphic structures on $E$ differing from 
a fixed $C^{\infty}$ one by an element of the Sobolev space 
$\Omega^{0,1}(\End)$. We shall denote by $E^{\del}$ the holomorphic bundle 
determined by any $\del\in\mathcal{A}$. Further, any such $\del$ also determines 
canonically a holomorphic structure on the bundle $E^*$, and hence a 
holomorphic structure on the bundle $\End \otimes \Omega^1_X$. By abuse of 
notation, we will denote all of these by $\del$. Given any Hermitian metric $H$
on $E$, we shall denote by $F_{\del, H}$ the curvature of the metric connection
on $E^{\del}$.

We shall denote by $\mathcal{G}^{\C}$ the complex gauge group of all complex 
automorphisms of $E$. So $\mathcal{G}^{\C}$ is the complexification of the 
unitary gauge group $\mathcal{G}$ of automorphisms that preserve a fixed 
Hermitian metric on $E$. The space $\mathcal{A}$ has a natural action of 
$\mathcal{G}^{\C}$ on it given by:
$$ g \cdot \del \quad = \quad g \circ\del\circ g^{-1} $$
for any $g \in \mathcal{G}^{\C}$.

Recall that a $GL_r(\C)$-Higgs bundle of rank $r$ and degree $d$ on $X$ 
is a pair of objects 
$(\del, \theta)$ where $\del \in \mathcal{A}$ as before is a holomorphic 
structure on $E$, while $\theta: 
E \rightarrow E \otimes \Omega^1_X$ is a holomorphic map (so $\theta \in 
H^0(\End^{\del} \otimes \Omega^1_X)$) 
such that $\theta \wedge \theta = 0$. In our case, since $X$ is a 
Riemann surface, the condition $\theta \wedge \theta = 0$ is vacuous. 
Since the focus of our attention will only be $GL_r(\C)$-Higgs bundles, we will 
drop the prefix and refer to these simply as Higgs bundles. Given a Higgs bundle 
$(\del, \theta)$, we say that a holomorphic subbundle $F \subset E^{\del}$ is 
$\theta$-invariant if $\theta(F) \subseteq F\otimes \Omega^1_X$. The Higgs 
bundle $(\del, \theta)$ is said to be semistable if
$$
\mu(E') \le \mu(E^{\del}) \quad \textrm{for all }\theta\textrm{-invariant 
proper holomorphic subbundles }E' \subset E^{\del},
$$
and it is said to be stable if the above inequality is strict.
As in the case of vector bundles, the complex gauge group $\mathcal{G}^{\C}$ 
acts naturally on the space of Higgs bundles over $X$:
$$ g \cdot (\del, \theta) \quad = \quad (g \circ\del\circ g^{-1}, g 
\circ\theta\circ g^{-1}) $$
for any $g \in \mathcal{G}^{\C}$. This action preserves the subset of stable 
Higgs bundles, and the quotient of this subset by the gauge group is the moduli 
space $\mathcal{H}(r,d)$ of Higgs bundles of rank $r$ and degree $d$ over $X$.

A stable pair of rank $r$ and degree $d$ on $X$ is a pair of objects $(\del, 
\phi)$ where $\del \in \mathcal{A}$ as before is a holomorphic structure on $E$, 
while $\phi$ is a holomorphic section of $E$ (so $\phi \in H^0(E^{\del})$). 
These were first studied by Bradlow in \cite{Br1} and \cite{Br2}.
Given any real number $\tau$, a pair $(\del, \phi)$ is said to be 
$\tau$-semistable if
\begin{align*}
\mu(E') &\le \tau \quad \textrm{for all holomorphic subbundles }E' \subseteq 
E^{\del}, \textrm{ and} \\
\mu(E^{\del}/E') &\ge \tau \quad \textrm{for all proper holomorphic subbundles 
}E' \subset E^{\del} \textrm{ that contain }\phi \\
&\qquad\quad\textrm { i.e., }\phi \in H^0(E').
\end{align*}
As usual, the pair is said to be $\tau$-stable if the above inequalities are 
strict. It can be shown that the set of $\tau$-semistable pairs is non-empty 
precisely when $\tau \in \left[\frac{d}{r}, \frac{d}{r-1}\right]$. For any fixed 
$\tau \in 
\left[\frac{d}{r}, \frac{d}{r-1}\right]$, the space of $\tau$-stable pairs over 
$X$ also 
has a natural action of the complex gauge group on it given by 
$$ g \cdot (\del, \phi) \quad = \quad (g \circ\del\circ g^{-1}, g \circ\phi ) $$
for any $g \in \mathcal{G}^{\C}$, and the resulting quotient is the moduli space 
of $\tau$-stable pairs over $X$.

Higgs bundles and stable pairs are both examples of augmented bundles that satisfy
the Hitchin-Kobayashi correspondence. In the case of Higgs bundles (proved by Hitchin)
this states that if a Higgs bundle $(\del, \theta)$ is stable, then the equation
\begin{equation*}
F_{\del, H} + [\theta, \theta^{*_H}] = \frac{d}{2r} i \omega {\bf I}
\end{equation*}
considered as an equation of the Hermitian metric $H$ on $E$ has a unique (up to 
scalar multiplication) smooth solution. Here $\theta^{*_H}$ is the adjoint of 
$\theta$ with respect to the metric $H$ defined by
$$ (\theta u, v)_H = (u, \theta^{*_H} v)_H, $$
and $[\cdot,\cdot]$ is the Lie bracket. 
Conversely, the existence of such a solution implies the polystability of the 
Higgs bundle. In the case of stable pairs, we have that if a pair $(\del, \phi)$
is $\tau$-stable, then 
\begin{equation*}
i\Lambda F_{\del, H} + \frac{1}{2}(\phi\otimes\phi^{*_H}) = \frac{\tau}{2}{\bf 
I}
\end{equation*}
considered as an equation of the Hermitian metric $H$ on $E$ has a unique (up to 
scalar multiplication) smooth solution. Here $\Lambda$ is the adjoint (with 
respect to 
the metric $H$) of $L$ (the Lefschetz operator given by the K\"ahler form 
$\omega$ on $X$). Again, there exists a suitable converse of this result. 

The Hitchin-Kobayashi correspondence has been generalized adequately for our 
purposes by Bradlow, Garcia-Prada, and Riera in \cite{BGI}. It is this interplay 
between stability and Hermitian metrics that allows us to interpret the various 
moduli spaces in more than one way. 

\section{Higgs Triples}
\setcounter{theorem}{0}
\setcounter{equation}{0}

We now introduce a new augmented bundle in gauge theory, which we shall use to 
study the moduli space of Higgs bundles. As before, given a Riemann surface $X$ 
and the complex vector bundle $E$ on $X$, we make the following definition.

\begin{defn}
A Higgs triple of rank $r$ and degree $d$ on $X$ is a 3-tuple of objects $\tri$ 
where $\del\in\mathcal{A}$ is a holomorphic structure on $E$, $\theta: E 
\rightarrow E \otimes \Omega^1_X$ is a holomorphic map (so $\theta$ is just a 
Higgs field), and $\phi$ is a holomorphic global section of $E^{\del}$ such that 
$\theta(\phi) = 0$.
\end{defn}

Alternately, a Higgs triple may be thought of as a pair of objects $(D'', 
\phi)$, where $D'' = \del+\theta$ is a Higgs bundle, and $\phi$ is a 
`holomorphic section of the Higgs bundle' i.e., $D''(\phi) = 0$. This viewpoint 
makes the analogy between stable pairs and these triples more obvious in the way 
they relate to vector bundles and Higgs bundles respectively.

In order to construct a moduli space, we define the space of $\T$ of Higgs 
triples:
\begin{multline*}
\T = \Big\{ \tri\in\mathcal{A} \times \Omega^0(\End \otimes \Omega^1_X) 
\times\Omega^0(E) 
\,\Big|\, \\
\del(\theta) = 0, \del(\phi) = 0 \textrm{ and } \theta(\phi) = 0 \Big\}. 
\end{multline*}
The condition $\theta(\phi) = 0$ implies that the line subbundle 
$[\phi] \subset E^{\del}$ generated by $\phi$ is $\theta$-invariant.

The space $\T$ admits a natural action of the complex gauge group 
$\mathcal{G}^{\C}$. This allows us to study these triples up to isomorphism. The 
action is given by
$$ g\cdot\tri \quad = \quad (g\circ \del \circ g^{-1}, g\circ \theta \circ 
g^{-1}, g\circ \phi) $$
for any $g \in \mathcal{G}^{\C}$.

Next, we define the notion of stability for these objects. Although this notion 
is understood more easily using a real parameter $\tau$, it may be defined 
intrinsically. In order to do this, we need the following. Let
\begin{align*}
\mu_M(\del, \theta) &= \textrm{Sup} \Big\{ \mu(E')\,\Big|\, E'\subseteq E^{\del} 
\textrm{ is a } \theta\textrm{-invariant holomorphic subbundle} \Big\}, \\
\mu_m\tri &= \textrm{Inf} \Big\{ \mu(E/E')\,\Big|\, E' \subset
E^{\del} \textrm{ is a } \theta\textrm{-invariant proper} \\
&\qquad\qquad\textrm{holomorphic subbundle containing } \phi \textrm{ i.e., } 
\phi\in 
\Omega^0(E') \Big\}.
\end{align*}
Now we can state the definition of stability of these triples.

\begin{defn}
A Higgs triple $\tri\in\T$ is said to be semistable if
$$ \mu_M(\del, \theta) \le \mu_m\tri. $$
The triple is said to be stable if the above inequality is strict. 
Moreover, 
if $\tau$ is any real number such that
$$ \mu_M(\del, \theta) \le \tau \le \mu_m\tri $$
then the triple $\tri$ will be said to be $\tau$-semistable. As before, the 
triple is said to be $\tau$-stable if the inequalities are strict.
\end{defn}

Clearly, a triple is semistable (resp. stable), if and only if there exists a 
$\tau$ such that it is $\tau$-semistable (resp. $\tau$-stable). Note that this 
notion of stability reduces to the notion of stability for stable pairs if the 
Higgs field vanishes ($\theta = 0$), and to the notion of stability for ordinary 
vector bundles if both the Higgs field and the global section vanish ($\theta = 
0$ and $\phi = 0$).

Higgs triples, like other augmented bundles in gauge theory, also satisfy the
Hitchin-Kobayashi correspondence. Suppose $H$ is a Hermitian metric on $E$. 
Let $\tri\in \T$ be any triple. We consider the equation
\begin{equation}
i\Lambda \left(F_{\del, H} + [\theta, \theta^{*_H}]\right) + 
\frac{1}{2}(\phi\otimes\phi^{*_H}) = \frac{\tau}{2}{\bf I}.
\end{equation}
Here $\phi^{*_H}$ is the adjoint of $\phi$, computed with respect to the metric 
$H$, and $\Lambda$ is the adjoint (with respect to the metric $H$) of $L$ (the 
Lefschetz 
operator given by $\omega$ on $X$), so that $\Lambda F_{D'', H}$ is an element 
of $\Omega^{0}(X, \End)$. $\tau$ is a real number, ${\bf I}$ is the identity 
section of End$(E)$, and all the terms of the equation take values in 
$\Omega^0(\End)$. Following the terminology from Bradlow \cite{Br1}, we will 
call this the $\tau$-vortex equation.

Note that in the absence of the Higgs field, the equation reduces to the vortex 
equation of Bradlow, while in the absence of the global section, the equation 
reduces to Hermitian-Yang-Mills equation of Hitchin, Simpson et al. As in these 
cases, the equation helps us determine certain preferred metrics on $E$. 

\begin{theorem}
\label{stability-vortex eqn}
Given $E$ and $X$ as above, let $\tri\in\T$ be any Higgs triple. Suppose that 
there exists a real number $\tau$ such that $\tri$ is a $\tau$-stable Higgs 
triple (as defined above). Then the $\tau$-vortex equation
$$ i\Lambda \left(F_{\del, H} + [\theta, \theta^{*_H}]\right) + 
\frac{1}{2}(\phi\otimes\phi^{*_H}) = \frac{\tau}{2}{\bf I} $$
considered as an equation for the Hermitian metric $H$ on $E$ has a unique 
smooth solution.

Conversely, suppose that for a given $\tau\in\R$, $\tau > 0$, there exists a 
Hermitian metric $H$ on $E$ such that the $\tau$-vortex equation is satisfied by 
$\tri$. Then $E^{\del}$ splits holomorphically as $E^{\del} = E_{\phi} \oplus 
E_s$ where
\begin{enumerate}
\item[(a)] $E_{\phi}$ is $\theta$-invariant, and contains the section $\phi$ (it 
is understood that $E_{\phi}$ obtains its holomorphic structure from $E^{\del}$)
\item[(b)] $(E_\phi, \theta|_{E_\phi}, \phi)$ is a stable Higgs triple with 
$\mu_M(E_\phi, \theta|_{E_\phi}) < \tau <$
\linebreak
$\mu_m(E_\phi, \theta|_{E_\phi}, 
\phi)$
\item[(c)] $E_s$, if non-empty, is a direct sum of $\theta$-invariant subbundles 
$E_i$
\item[(d)] all the Higgs bundles $(E_i, \theta_i)$ are stable (in the usual 
sense of stability for Higgs bundles) with $\mu(E_i) = \tau$, where $\theta_i = 
\theta|_{E_i}$
\end{enumerate}
\end{theorem}

We point out that a split $E^{\del} = E_{\phi} \oplus E_s$ cannot occur unless 
$\mu(E_s) = \tau$ is a rational number with denominator less than rk$(E)$. 
Hence, for {\it generic} values of $\tau$, $E_s = \emptyset$.

Although the above theorem can be proved using methods similar to those of 
Simpson \cite{Si1} and Bradlow \cite{Br2} (done independently by the author, 
unpublished), it also follows as a special case from the recent work of 
Bradlow, Garcia-Prada and Riera \cite{BGI}. The essential observation is that 
the 
moment map 
for the action of the unitary gauge group $\mathcal{G}$ on the K\"ahler space
$\mathcal{A} \times \Omega^1(\End) \times \Omega^0(E)$ is
$$ \left(F_{\del, H} + [\theta, \theta^{*_H}]\right)
 - \frac{i}{2}(\phi\otimes  \phi^{*_H})\omega, $$
where $H$ is a fixed Hermitian metric on $E$.

The values of the parameter $\tau$ for which the solution set of the vortex 
equation is non-empty lie in a closed interval as in \cite{BDW}.

\begin{prop}
\label{imageoff}
There is a solution to the $\tau$-vortex equation only if $\D\frac{d}{r} \le 
\tau \le \frac{d}{r-1}$.
\end{prop}
\begin{proof}
By taking the trace of the $\tau$-vortex equation and integrating over the 
Riemann surface $X$, we get
$$ d + \frac{1}{2}\parallel\phi\parallel^2 = r\tau, $$
which implies that there is no solution unless
$$ \frac{d}{r} \le \tau. $$
For an upper bound, consider the $\theta$-invariant line subbundle $[\phi] 
\subset E^{\del}$. If the vortex equation is satisfied, then either $\tri$ is 
$\tau$-stable, or it splits holomorphically into a direct sum of $\tau$-stable 
pieces with prescribed slope. In the latter case, the $\tau$-vortex equation is 
then satisfied by each piece separately so that either way $\tri$ must be 
$\tau$-semistable. Thus we have
$$ \tau \le \frac{\textrm{deg}(E/[\phi])}{\textrm{rk}(E/[\phi])} $$
and since $\textrm{deg}([\phi]) \ge 0$, we conclude that
$$ \tau \le \frac{d}{r-1}. $$
\end{proof}

The converse to the above proposition is also true, and will be proved
later.

\section{Moduli Spaces}

From the previous section we see that for generic values of $\tau$, the set
$$ \mathcal{V}_{\tau} = \left\{ \tri \in \T \, \Big| \, \Lambda F_{D'',
H} - \frac{i}{2}(\phi\otimes\phi^{*_H}) = -i\frac{\tau}{2}{\bf I}
\textrm{ for some metric } H \right\} $$
consists exactly of the $\tau$-stable Higgs triples. Further, upon fixing a 
Hermitian metric $H$ on $E$, the unitary gauge group $\mathcal{G}$ acts 
symplectically on $\T$, and the moment map for this action is given by
$$ \Psi\tri = \Lambda F_{D'', H} - \frac{i}{2}(\phi\otimes\phi^{*_H}). $$
Hence for generic $\tau$, the spaces $\mathcal{V}_\tau$ consist precisely of all 
the $\mathcal{G}^\C$ orbits through the triples in $\Psi^{-1}(-
i\frac{\tau}{2}{\bf I})$. As a result, we may define the moduli space of $\tau$-
stable Higgs triples (for generic $\tau$) in two different ways:
$$ \B_\tau = \mathcal{V}_\tau/\mathcal{G}^\C = \Psi^{-1}(-i\frac{\tau}{2}{\bf 
I})/\mathcal{G}. $$

Along the lines of Bradlow et al. \cite{BDW}, we construct a `master space' 
which will contain stable and semi-stable Higgs triples for all values of 
$\tau$. To do this, as in \cite{BDW}, we use a 
moment map to construct $\hatb$: we replace the full unitary gauge group 
$\mathcal{G}$ by a subgroup $\mathcal{G}_0$ which has a $U(1)$ quotient and 
whose Lie algebra is the $L^2$ ortho-complement of the constant multiples of the 
identity. Denoting the new moment map by $\Psi_0$, we then obtain $\hatb$ 
as $\Psi_0^{-1}(0)/\mathcal{G}_0$ symplectically. Finally, we obtain the 
complex structure on $\hat{\mathcal{B}}$ by looking at $\mathcal{G}_0^\C$ orbits
through the triples in $\Psi_0^{-1}(0)$.
For details on the construction of the subgroups $\mathcal{G}_0$ and
$\mathcal{G}_0^\C$, we refer the reader to section 2.2 in \cite{BDW}.

The complex structure of $\hatb$ is essentially obtained from considering the 
action of $\mathcal{G}_0^\C$ on an appropriate subset of $\T$. Consequently we 
first focus our attention on the action of $\mathcal{G}^\C$ on this subset, a 
slight modification of which will lead us to our goal. First, we define
\begin{equation*}
\T^* = \left\{ \tri\in\T \,\Big|\, (\del, \theta) \textrm{ is semistable if }
\phi = 0 \right\}.
\end{equation*}
Then $\T^*$ is an open subset of $\T$. Next, for any $\del' \in \mathcal{A}$, 
let
$$ \T^*_{\del'} = \left\{ \tri\in\T^* \,|\, \del = \del' \right\}. $$
We consider the infinitesimal deformations of pairs $(\theta, \phi)$ inside the 
space $\T^*_{\del}$ for a fixed $\del \in \mathcal{A}$. In order to do this, we 
look at the complex
\begin{equation*}
(\End \otimes K) \oplus E \longrightarrow E \otimes K,
\end{equation*}
where the map is given by
$$ (v,w) \mapsto v\phi + \theta w. $$
Let ${\bf N}^{\tri}$ be the resulting long exact sequence. We define
\begin{equation*}
\T^{**}_{\del} = \left\{ \tri\in\T \,\Big|\, H^2({\bf N}^{\tri}) = 0 \right\}.
\end{equation*}
Then $\T^{**}_{\del}$ is an open subset of $\T^*_{\del}$. Now let
\begin{equation*}
\T^{**} = \bigcup\limits_{\del\in\mathcal{A}} \T^{**}_{\del}.
\end{equation*}
This space $\T^{**}$ is an open subset of $\T^*$. For purely technical reasons 
(which will become clear later), we shall restrict our attention to this subset 
of $\T$.

Now we use the standard approach of infinitesimal deformations 
to compute the obstruction to $\T^{**}/\mathcal{G}^{\C}$ being a manifold. These 
arise from the hypercohomology of the complex ${\bf End}^{\tri}$:
\begin{equation*}
\End \longrightarrow (\End \otimes K) \oplus E \longrightarrow E \otimes K,
\end{equation*}
where the maps are given by
\begin{equation*}
u \mapsto ([u, \theta], u\phi) \qquad \textrm{and} \qquad (v, w) \mapsto
v\phi + \theta w,
\end{equation*}
respectively. The complex which allows us to compute hypercohomology then 
becomes
\begin{align*}
0 \longrightarrow \Omega^0(\textrm{End }E) \overset{d_1}\longrightarrow
& \, \Omega^{0,1}(\textrm{End }E) \oplus \Omega^{1,0}(\textrm{End }E)
\oplus \Omega^0(E) \\
& \overset{d_2}\longrightarrow \Omega^{1,1}(\textrm{End }E) \oplus
\Omega^{0,1}(E) \oplus \Omega^{1,0}(E) \overset{d_3} \longrightarrow
\Omega^{1,1}(E) \longrightarrow 0,
\end{align*}
where the maps $d_1$, $d_2$ and $d_3$ are given by
\begin{align*}
d_1(u) &= (-\del(u), [u, \theta], u\phi) \\
d_2(\alpha, \beta, \gamma) &= ( [\alpha, \theta] + \del(\beta), \alpha\phi + 
\del(\gamma), \beta\phi + \theta\gamma ) \\
d_3(\lambda, \rho, \sigma) &= \lambda\phi + \theta\rho + \del(\sigma).
\end{align*}
We shall call the above complex ${\bf C}^{\tri}$.

\begin{prop}
The cup product from $H^1({\bf C}^{\tri}) \times H^1({\bf C}^{\tri})
\rightarrow$ \linebreak $H^2({\bf C}^{\tri})$ vanishes.
\end{prop}

\begin{proof}
This follows now from smoothness of the Higgs moduli space and 
the moduli space of stable pairs, and our construction of $\T^{**}$.
\end{proof}

Now we look at the action of $\mathcal{G}_0^\C$ on $\T^{**}$. This requires us 
to restrict to the following subcomplex of ${\bf C}^{\tri}$:
\begin{align*}
0 \longrightarrow \Omega^0(\textrm{End }E)_0 &\overset{d_1} \longrightarrow \, 
\Omega^{0,1}(\textrm{End }E) \oplus \Omega^{1,0}(\textrm{End }E) \oplus 
\Omega^0(E) \\
& \overset{d_2}\longrightarrow \Omega^{1,1}(\textrm{End }E) \oplus 
\Omega^{0,1}(E) \oplus \Omega^{1,0}(E) \overset{d_3} \longrightarrow
\Omega^{1,1}(E) \longrightarrow 0,
\end{align*}
which we shall denote by ${\bf C}^{\tri}_0$. The only difference is in the first 
term which is now $\Omega^0(\End)_0$. Recall that $\Omega^0(\End)_0$ was the 
$L^2$-orthogonal complement of the constant multiples of the identity in 
$\Omega^0(\End)$. The important point to note here is that the cup 
product from $H^1({\bf C}^{\tri}_0) \times H^1({\bf C}^{\tri}_0) 
\rightarrow H^2({\bf C}^{\tri}_0)$ still vanishes. Consequently we have the 
following.

\begin{prop}
$\T^{**}$ is a smooth submanifold of $\mathcal{A} \times \Omega^{0,1}(E) \times 
\Omega^0(E)$.
\end{prop}

\begin{defn}
We define a triple $\tri \in \T^{**}$ to be simple if $H^0(C^{\tri}_0) =
0$. Denote by $\T_{\sigma}$ the set of all simple triples in $\T^{**}$.
\end{defn}
Note that $\T_\sigma$ is open in $\T^{**}$. This allows us to conclude the 
following.
\begin{prop}
\label{connectone}
$\T_\sigma/\mathcal{G}_0^\C$ is a complex manifold. Moreover, we have the 
identification 
$$ T_{[\del, \theta, \phi]}(\T_\sigma/\mathcal{G}_0^\C) = H^1({\bf 
C}_0^{\tri}). $$
\end{prop}

We now study the symplectic structure on $\T_\sigma/\mathcal{G}_0^\C$.

\begin{prop}
The moment map for the action of the subgroup $\mathcal{G}_0$ on $\T^{**}$ is 
given by
\begin{equation*}
\Psi_0\tri = \pi^{\perp}\Psi\tri = \Psi\tri - \frac{1}{r}\int_X \textrm{Tr } 
\Psi\tri\cdot {\bf I}
\end{equation*}
\end{prop}
\begin{proof}
We observe that $\Psi_0 = j^*\Psi$ where $j: \mathcal{G}_0 \rightarrow 
\mathcal{G}_0$ is the inclusion. For details, see section 2.4 in \cite{BDW}.
\end{proof}

We are now ready to define the master space $\hatb$.

\begin{defn}
Let
\begin{equation*}
\hatb = (\Psi_0^{-1}(0) \cap \T^*)/\mathcal{G}_0
\end{equation*}
be the Marsden-Weinstein reduction by the symplectic action of $\mathcal{G}_0$ 
and
\begin{equation*}
\hatb_0 = (\Psi_0^{-1}(0)\cap \T_\sigma)/\mathcal{G}_0.
\end{equation*}
\end{defn}

\begin{prop}
\label{connecttwo}
The space $\hatb$ is a Hausdorff topological space. The space $\hatb_0$ 
is a Hausdorff symplectic manifold.
\end{prop}
\begin{proof}
The first follows in an identical manner to the similar result in \cite{BDW}, 
while the second is a consequence of the Marsden-Weinstein reduction theorem for 
Banach spaces.
\end{proof}

Finally, in order to obtain the complex structure on $\hatb_0$, we define
\begin{defn}
Let $\mathcal{V}_0 \subset \T^{**}$ denote the subset of $\mathcal{G}_0^{\C}$ 
orbits through points in $\Psi_0^{-1}(0)$, that is
$$ \mathcal{V}_0 = \left\{ \tri\in\T^{**} \, \Big|\, \Psi_0(g\tri) = 0 \textrm{ 
for some } g\in\mathcal{G}_0^{\C} \right\}. $$
\end{defn}

Note that $\mathcal{V}_0 \cap \T_{\sigma}$ is an open subset of $\T_{\sigma}$. 
Using the same technique as in \cite{BDW}, it also follows that $\mathcal{V}_0 
\cap \T_{\sigma}$ is connected. Further, there is a bijective correspondence 
between $\hatb_0$ and $(\mathcal{V}_0 \cap \T_{\sigma})/\mathcal{G}_0^{\C}$, so 
that using propositions \ref{connectone} and \ref{connecttwo}, we obtain the 
following.

\begin{prop}
$\hatb_0 = (\mathcal{V}_0 \cap \T_{\sigma})/\mathcal{G}_0^{\C}$ is a smooth, 
Hausdorff manifold.
\end{prop}

\section{$S^1$ Action and Morse Theory}
\setcounter{theorem}{0}
\setcounter{equation}{0}

A key feature of the master space $\hatb$ is that it carries a natural 
$S^1$-action. This arises from the quotient $\mathcal{G}/\mathcal{G}_0 = U(1)$.
The action is given by:
\begin{equation*}
e^{i\rho}\cdot [\del, \theta, \phi] = [\del, \theta, g_\rho \phi].
\end{equation*}
Here $g_\rho$ denotes the gauge transformation diag$(e^{i\rho/r}, \ldots, 
e^{i\rho/r})$. As shown in \cite{BDW}, this action is well-defined and 
independent of the choice of the $r$-th root of unity, for if $h = e^{2\pi i/r} 
\cdot {\bf I}$, then $h \in \mathcal{G}_0$ and
$$ [\del, \theta, h\phi] = [h^{-1}\del, h^{-1}\theta, \phi] = [\del, \theta, 
\phi]. $$

\begin{prop}
The action of $U(1)$ on $\hatb_0$ is holomorphic and symplectic. The moment map 
for the action is given by:
\begin{equation*}
\hat{f}[\del, \theta, \phi] = -2\pi i \Big(
\frac{\parallel\phi\parallel^2}{4\pi r} + \mu(E) \Big).
\end{equation*}
Further, this action extends continuously to $\hat{\mathcal{B}}$ as does the 
moment map $\hat{f}$.
\end{prop}
\begin{proof}
The proof is identical to the one in \cite{BDW}. The only observation which one 
needs is that the trace of the term contributed by the Higgs field is zero i.e.,
$$ \textrm{Tr}([\theta, \theta^{*_H}]) = 0. $$
\end{proof}

For convenience, we define $f:\hatb \rightarrow \R$ by
\begin{equation*}
f = -\frac{1}{2\pi i} \hat{f}.
\end{equation*}
The moment map $\hat{f}$ is the same as the one obtained in \cite{BDW}. The next 
proposition summarizes the essential properties of $f$. The proofs are identical 
to proposition (2.14) in \cite{BDW}.

\begin{prop}
\begin{enumerate}
\item[(i)] The image of $f$ is the interval $\left[\frac{d}{r}, 
\frac{d}{r-1}\right]$.
\item[(ii)] The critical points of $f$ on $\hatb_0$ are exactly the fixed points 
of the $U(1)$ action. The critical values of $f$ are precisely the image under 
$f$ of the fixed point set of the $U(1)$ action on $\hatb$.
\item[(iii)] If $\tau$ is a regular value of $f$, then the space 
$f^{-1}(\tau)/U(1)$ is $\B_\tau$, the moduli space of $\tau$-stable Higgs 
triples.
\end{enumerate}
\end{prop}
\begin{proof}
(i) $\tau$ is in the image of $f$ if and only if the equation $\Psi(\tri) = 
-\frac{i \tau}{2}{\bf I}$ has a solution. By proposition \ref{imageoff}, the 
range for $\tau$ is in $[d/r, d/(r-1)]$. Next, proposition \ref{endpoints} gives 
explicit elements of the spaces at the end points, $\B_{\frac{d}{r}}$ and 
$\B_{\frac{d}{r-1}}$. Therefore, since $\mathcal{V}_0 \cup \T_{\sigma}$ is 
connected, the result follows. (ii) This follows from the fact that $\hat{f}$ is 
a moment map for the $S^1$ action. (iii) See \cite{BDW}.
\end{proof}

Next, we look at the level sets $f^{-1}(\tau)$ when $\tau$ 
is a critical value in $\left[\frac{d}{r}, \frac{d}{r-1}\right]$. We 
start with a simple but important observation.

\begin{prop}
\label{endpoints}
The level set corresponding to the minimum is precisely the moduli space
of semistable Higgs bundles of degree $d$ and rank $r$:
$$ f^{-1}\left(\frac{d}{r}\right) = \mathcal{M}_{Higgs}(r,d), $$
while the level set corresponding to the maximum is precisely the moduli
space of semistable Higgs bundles of degree $d$ and rank $r-1$:
$$ f^{-1}\left(\frac{d}{r-1}\right) = \mathcal{M}_{Higgs}(r-1,d). $$
\end{prop}
\begin{proof}
By construction, $f^{-1}(d/r)$ consists exactly of $(d/r)$-semistable triples. 
However, the standard trick of taking the trace of the vortex equation and 
integrating over $X$ shows that $\phi$ is forced to be 
identically $0$ when $\tau = d/r$. Then $(d/r)$-semistable is the same as 
semistability in the usual sense of a Higgs bundle. To see this for the maximum 
level set, consider a solution $\tri$ of the vortex equation with $\tau = d/(r-
1)$. Then by theorem \ref{stability-vortex eqn}, $\tri$ is either stable or 
splits holomorphically. If the triple is stable, then since the line subbundle
$[\phi] \subset E^{\del}$ is $\theta$-invariant, the stability criteria gives us
$$ \mu(E/[\phi]) > d/(r-1) \quad \Rightarrow \quad \textrm{deg}([\phi]) < 0, $$
which is a contradiction. Hence the triple $(D'', \phi)$ splits. But we also 
have
$$ \mu(E/[\phi]) = d/(r-1) \quad \Rightarrow \quad \textrm{deg}([\phi]) = 0. $$
Theorem \ref{stability-vortex eqn} and its proof in the converse direction now 
forces $\phi$ to be a constant section of a trivial line subbundle. This means 
that $E^{\del}$ splits as $E^{\del} = \mathcal{O} \oplus E_s$, where $(E_s, 
\theta|_{E_s})$ is a semistable Higgs bundle of degree $d$ and rank $r-1$. 
Further, as in \cite{BDW}, due to theorem \ref{stability-vortex eqn} the Higgs 
bundle $(E_s, \theta|_{E_s})$ is a direct sum of stable Higgs bundles all of the 
same slope, so that it is equal to its graded $S$-equivalence class in 
$\mathcal{M}_{Higgs}(r-1, d)$. Thus the map $\tri = (\mathcal{O} \oplus E_s 
\mapsto (E_s, \theta|_{E_s})$ establishes the correspondence between $f^{-
1}(d/(r-1))$ and $\mathcal{M}_{Higgs}(r-1, d)$.
\end{proof}

We now turn our attention to the level sets corresponding to the critical values 
in the interior of the interval $\left[\frac{d}{r}, \frac{d}{r-1}\right]$. To 
this end we 
make the following definition.

\begin{defn}
Let $Fix(\hatb)$ be the $U(1)$ fixed point set in $\hatb$. Then for any critical 
value $\tau \in \left[\frac{d}{r}, \frac{d}{r-1}\right]$, define
$$ \mathcal{Z}_\tau = f^{-1}(\tau) \cap Fix(\hat{\mathcal{B}}). $$
\end{defn}

Let $\tau = \frac{p}{q}$ be a critical value such that $\frac{d}{r} < \tau < 
\frac{d}{r-1}$. Then for any $\tri \in \mathcal{Z}_\tau$, we have $\mu_M(\del, 
\theta) = \frac{p}{q} = \mu_m(\tri)$. Further, theorem \ref{stability-vortex 
eqn} allows us to make certain immediate observations which we summarize in the 
next proposition.

\begin{prop}
\label{tuples lemma}
Let $\frac{p}{q} \in (\frac{d}{r}, \frac{d}{r-1})$ be a critical value, and 
$\tri \in f^{-1}(\frac{p}{q})$. As in theorem \ref{stability-vortex eqn}, 
since $\tri$ satisfies the vortex equation, we know that $E^{\del}$ splits 
holomorphically as $E^{\del} = E_\phi \oplus E_s$ where $E_s$ is possibly a 
direct sum of bundles $E_i$. Let $(r_\phi, d_\phi)$ and $(r_i, d_i)$ be the rank 
and degree of $E_\phi$ and $E_i$, respectively. Suppose $\tri$ is not 
$\frac{p}{q}$-stable. Then the following holds:
\begin{enumerate}
\item[(i)] $\D \frac{d_i}{r_i} = \frac{d - d_\phi}{r - r_\phi} = \frac{p}{q}$
\item[(ii)] $r_\phi + \sum\limits_{i} r_i = r$
\item[(iii)] $\D \frac{d_\phi}{r_\phi} < \frac{p}{q} < \frac{d_\phi}{r_\phi - 
1}$
\end{enumerate} 
Conversely, given any stable Higgs triple $(E_\phi, \theta_\phi, \phi)$  and 
stable Higgs bundles $(E_i, \theta_i)$ so that the above conditions are 
satisfied, we have a representative for a fixed point $(\del, \theta, 
\phi)$ in $f^{-1}(\frac{p}{q})$ where $(E^{\del}, \theta) = (E_\phi, 
\theta_\phi) \oplus \bigoplus\limits_i (E_i, \theta_i)$.
\end{prop}
\begin{proof}
In the forward direction, since $\tri$ is not $\frac{p}{q}$-stable but just 
semistable, (i) and (ii) are obvious consequences of theorem
\ref{stability-vortex eqn}. (iii) follows from applying the discussion on the 
possible range of the parameter $\tau$ to the $\frac{p}{q}$-stable triple 
$(\bar{\partial}_{E_\phi}, \theta_\phi, \phi)$.

To see the converse, we use the same key theorem to find metrics on the various 
pieces $E_\phi$ and the bundles $E_i$ so that the $\frac{p}{q}$ vortex equation 
holds for the various pieces. Then we piece these metrics together to obtain a 
metric $H$ on $E^{\del} = E_\phi \oplus \bigoplus\limits_i E_i$ so that the 
resulting triple $\tri$ satisfies the $\frac{p}{q}$-vortex equation. However, 
since this triple is only $\frac{p}{q}$-semistable, it is in $f^{-
1}(\frac{p}{q})$.
\end{proof}

In order to investigate the level sets at the critical points, we stratify them 
in the following manner. Given any critical value $\tau = \frac{p}{q}$, let
\begin{multline*}
\mathcal{I}_\tau = \Big\{(d_\phi, r_\phi, r_1, \ldots, r_n)\in \mathbb{Z}^{n+2} 
\,\Big| \, \frac{p}{q} = \tau = \frac{d-d_\phi}{r-r_\phi} \textrm{ and 
conditions (ii), (iii)} \\
\textrm{of the previous proposition are satisfied} \Big\}.
\end{multline*}

Note that given such an $(n+2)$-tuple $(d_\phi, r_\phi, r_1, \ldots, r_n)$, the 
degrees $d_i$ are determined uniquely using condition (i) in the previous 
proposition.

This allows us to write $\mathcal{Z}_\tau$ as a disjoint union of sets 
$\mathcal{Z}(d_\phi, r_\phi, r_1, \ldots, r_n)$ where
$$ \mathcal{Z}(d_\phi, r_\phi, r_1, \ldots, r_n) = \left\{ \tri 
\in \mathcal{Z}_\tau \, \Big| \, (E^{\del}, \theta) = (E_\phi, \theta_\phi) 
\oplus \bigoplus\limits_i (E_i, \theta_i) \right\}. $$

Hence we get
\begin{equation*}
\mathcal{Z}_\tau = \bigcup\limits_{(d_\phi, r_\phi, r_1, \ldots, r_n) \in 
\mathcal{I}_\tau} \mathcal{Z}(d_\phi, r_\phi, r_1, \ldots, r_n).
\end{equation*}

By using the following convention we can include in the above discussion the 
critical values at the end points of the interval i.e., $\tau = \frac{d}{r}$ and 
$\tau = \frac{d}{r-1}$. When $\tau = \frac{d}{r}$, we will set $d_\phi = 0$ and 
$r_\phi = 0$, and when $\tau = \frac{d}{r-1}$, we will set $d_\phi = 0$ and 
$r_\phi = 1$.

\section{Algebraic Stratification}
\setcounter{theorem}{0}
\setcounter{equation}{0}

In this section we define a natural stratification of $\hat{\mathcal{B}}$. We do 
this by defining two different filtrations on any Higgs triple in the master 
space, which arise as a consequence of the definition of stability for such 
objects. These filtrations are similar to the ones defined in \cite{BDW}, except 
we require that the relevant subbundles be invariant under the Higgs field. 

\begin{prop}
\label{minus filtration}
(The $\mu_-$-filtration). Let $\tri$ be a stable Higgs triple. There is a 
filtration of $E$ by 
$\theta$-invariant subbundles 
\begin{equation}
0 \subset E_\phi = F_0 \subset F_1 \subset \cdots \subset F_n = E
\end{equation}
such that the following conditions hold:
\begin{enumerate}
\item[(a)] $\phi \in H^0(E_\phi)$, $(E_\phi, \theta_\phi, \phi)$ is a stable 
Higgs triple (where $\theta_\phi = \theta|_{E_\phi}$), and
$\mu_M(E_\phi, 
\theta_\phi) < \mu_m(E, \theta, \phi) < \mu_m(E_\phi, \theta_\phi, \phi)$,
\item[(b)] for $1 \le i \le n$, the quotients $(F_i/F_{i-1}, \theta_i/ 
\theta_{i-1})$ are stable Higgs bundles each of slope $\mu_m(E, \theta, 
\phi)$ (where $\theta_i = \theta_{F_i}$),
\item[(c)] $E_\phi$ has minimal rank such that (a) and (b) are satisfied.
\end{enumerate}
Consequently, the subbundle $E_\phi$ is uniquely determined and the graded 
object
$$ gr^-(\del, \theta) = (E_\phi, \theta_\phi) \oplus (F_1/F_0, 
\theta_1/\theta_0) \oplus \cdots \oplus (F_n/F_{n-1}, \theta_n/\theta_{n-1}) $$
is unique up to isomorphism of $(F_1/F_0, \theta_1/\theta_0) \oplus \cdots 
\oplus 
(F_n/F_{n-1}, \theta_n/\theta_{n-1})$.
\end{prop}

Before we can prove this, we need the following.

\begin{lemma}
Let $\tri$ be a stable triple. Let $E_\phi \subset E^{\del}$ be a 
$\theta$-invariant subbundle such that $\phi \in H^0(E_\phi)$ and $\mu(E/E_\phi) 
= \mu_m(E, \theta, \phi)$. Let $\theta_\phi = \theta|_{E_\phi}$. Then
\begin{enumerate}
\item[(a)] $\mu_M(E_\phi, \theta_\phi) \le \mu_M(E, \theta)$,
\item[(b)] $\mu_m(E_\phi, \theta_\phi, \phi) \le \mu_m(E, \theta, \phi)$ and 
the inequality is strict if $E_\phi$ has minimal rank among all subbundles 
satisfying the hypotheses of this Lemma,
\item[(c)] the Higgs triple $(E_\phi, \theta_\phi, \phi)$ is stable,
\item[(d)] $(E/E_\phi, \theta/\theta_\phi)$ is a semistable Higgs bundle,
\item[(e)] $\mu(E_\phi) < \mu_m(E, \theta, \phi)$,
\item[(f)] If $E_\phi$ has minimal rank among all the $\theta$-invariant 
subbundles satisfying the hypotheses of this Lemma, and $E'_\phi$ is any other 
$\theta$-invariant subbundle satisfying the same, then $E_\phi \subset E'_\phi$.
\end{enumerate}
\end{lemma}

\begin{proof}
\begin{enumerate}
\item[(a)] This follows by definition of $\mu_M$.
\item[(b)] The proof here is identical to \cite{BDW} with a small modification: 
we consider only Higgs field invariant subbundles. This however, does not affect 
the argument.
\item[(c)] Using (a), (b) and the fact that $(E, \theta, \phi)$ is stable, it 
follows that $(E_\phi, \theta_\phi, \phi)$ is stable.
\item[(d)] Suppose $(E/E_\phi, \theta/\theta_\phi)$ is not semistable as a Higgs 
bundle. Then there is some $(\theta/\theta_\phi)$-inv. subbundle $F \subset 
E/E_\phi$ such that $\mu(F) = \mu_M(E/E_\phi, \theta/\theta_\phi)$
\linebreak 
$> \mu(E/E_\phi)$. 
Consider the Higgs extension of $(E_\phi, \theta_\phi)$ by $(F, \theta_F)$ where 
$\theta_F = (\theta/\theta_\phi)|_F$. This gives us an exact sequence of Higgs 
bundles
$$ 0 \longrightarrow (E_\phi, \theta_\phi) \longrightarrow (E', \theta') 
\longrightarrow (F, \theta_F) \longrightarrow 0. $$
Recall that an exact sequence of Higgs bundles as above is simply two exact 
sequences with maps between them so that the resulting diagram commutes:
\begin{equation*}
\begin{CD}
0 @>>> E_\phi @>>> E' @>>> F @>>> 0 \\
&& @VV\theta_\phi V @VV\theta' V @VV\theta_F V \\
0 @>>> E_\phi \otimes K @>>> E' \otimes K @>>> F \otimes K @>>> 0\\
\end{CD}
\end{equation*}
where $K$ is the canonical bundle on the Riemann surface $X$.
Note that by construction $E'$ is $\theta$ invariant, and we let $\theta' = 
\theta|_{E'}$. The rest of the argument works the same way as in \cite{BDW}.
\item[(e)] Since $(E, \theta, \phi)$ is stable, $\mu(E_\phi) \le \mu_M(E, 
\theta) < \mu_m(E, \theta, \phi)$.
\item[(f)] Suppose $E_\phi$ and $E'_\phi$ are as in statement (f) of the 
proposition. Let $\theta'_\phi = \theta|_{E'_\phi}$. We observe that the 
inclusion $(E_\phi, \theta_\phi) \rightarrow (E, \theta)$ and projection $(E, 
\theta) \rightarrow (E/E'_\phi, \theta/\theta'_\phi)$ are both morphisms of 
Higgs bundles. Hence, so is the composition $(E_\phi, \theta_\phi) 
\longrightarrow (E/E'_\phi, \theta/\theta'_\phi)$. Taking the kernel and the 
image of this map we get the following exact sequence of Higgs bundles:
$$ 0 \longrightarrow (N, \theta_N) \longrightarrow (E_\phi, \theta_\phi) 
\longrightarrow (L, \theta_L) \longrightarrow 0, $$
where $N$ and $L$ are respectively the kernel and the image, $\theta_N = 
\theta_\phi|_N$ and $\theta_L = (\theta/\theta'_\phi)|_L$. The rest of 
the argument proceeds as in \cite{BDW}.
\end{enumerate}
\end{proof}

\begin{proof}[Proof of Proposition \ref{minus filtration}]
The previous Lemma allows us to find a unique $\theta$-invariant subbundle 
$E_\phi \subset E$ of minimal rank such that
\begin{enumerate}
\item[(i)] $\phi \in H^0(E_\phi)$,
\item[(ii)] $\mu_M(E_\phi, \theta_\phi) < \mu_m(E, \theta, \phi) < \mu_m(E_\phi, 
\theta_\phi, \phi)$ (in particular $(E_\phi, \theta_\phi, \phi)$ is a stable 
Higgs triple),
\item[(iii)] $\mu(E/E_\phi) = \mu_m(E, \theta, \phi)$,
\item[(iv)] $(E/E_\phi, \theta/\theta_\phi)$ is a semistable Higgs bundle.
\end{enumerate}
Now let
$$ 0 \subset (Q_1, \eta_1) \subset \cdots \subset (Q_n, \eta_n) = (E/E_\phi, 
\theta/\theta_\phi) $$
be the Harder-Narasimhan (HN) filtration for the Higgs bundle $(E/E_\phi, 
\theta/\theta_\phi)$. Note 
that the bundles $Q_i$ are all $\theta/\theta_\phi$-invariant and $\eta_i = 
(\theta/\theta_\phi)|_{Q_i}$. Let $\pi: (E, \theta) \rightarrow (E/E_\phi, 
\theta/\theta_\phi)$ be the projection map (which as discussed earlier is a 
morphism of Higgs bundles). Define $(F_i, \theta_i) = \pi^{-1}(Q_i, \eta_i)$. 
This gives us the filtration we seek.
\end{proof}

This allows us to make the following definition.

\begin{defn}
For any stable Higgs triple $\tri$, we define the $\mu_-$ grading to be given by
$$ gr^-\tri = (gr^-(\del, \theta), \phi), $$
where $gr^-(\del, \theta)$ is as above.
\end{defn}

\begin{prop}
(The $\mu_+$-filtration) Let $\tri$ be a stable Higgs triple. There is a 
filtration of $E$ by $\theta$-invariant subbundles 
\begin{equation}
0 = F_0 \subset F_1 \subset \cdots \subset F_n \subset F_{n+1} = E
\end{equation}
such that the following conditions hold: if $(E, \theta)$ is a semistable Higgs 
bundle, then this is the usual generalization of the HN filtration to Higgs 
bundles so that the bundles $(F_i/F_{i-1}, \theta_i/\theta_{i-1})$ are stable 
Higgs bundles of slope $= \mu_M(E, \theta) = \mu(E)$. Otherwise
\begin{enumerate}
\item[(a)] for $1 \le i \le n$, the quotients $(F_i/F_{i-1}, 
\theta_i/\theta_{i-1})$ are stable Higgs bundles each of slope $\mu_M(E, 
\theta)$,
\item[(b)] $\phi$ projects to some non-zero $\psi \in H^0(E/F_n)$ and 
$(E/F_n, \theta/\theta_n, \psi)$ is a stable Higgs triple, and $\mu_M(E/F_n, 
\theta/\theta_n) < \mu_M(E, \theta) < \mu_m(E/F_n, \theta/\theta_n, \psi)$,
\item[(c)] $E/F_n$ has minimal rank such that (a) and (b) are satisfied.
\end{enumerate}
When $(E, \theta)$ is unstable, the subbundle $E/F_n$ is uniquely determined by 
the graded object
$$ gr^+(\del, \theta) = (E/F_n, \theta/\theta_n) \oplus (F_1/F_0, \theta_1/
\theta_0) \oplus \cdots \oplus (F_n/F_{n-1}, \theta_n/\theta_{n-1}) $$
and is unique up to isomorphism of $(F_1/F_0, \theta_1/\theta_0) \oplus \cdots 
\oplus (F_n/F_{n-1}, \theta_n/\theta_{n-1})$.
\end{prop}

\begin{proof}
If $(\del, \theta)$ is a semistable Higgs bundle, then $\mu_M(\del, \theta) = 
\mu(E)$ and we simply use the usual HN filtration of the Higgs bundle $(\del, 
\theta)$. If $(\del, \theta)$ is unstable, then let $F$ be the unique maximal 
destabilizing subbundle of $E^{\del}$. Note that $F$ is $\theta$-invariant and 
$(F, \theta|_F)$ is a semistable Higgs bundle. In this case we let $0 = F_0 
\subset \cdots \subset F_n = F$ be the usual HN filtration for $F$ and thus 
obtain the filtration as stated in the proposition.

Since $\mu(F) = \mu_M(\del, \theta)$ by the choice of $F$, part (a) follows 
directly. Part (b) follows from a similar argument to the corresponding 
statement in \cite{BDW}. The only difference is that we consider Higgs field 
invariant subbundles instead of any subbundles, and all extensions are 
extensions of Higgs bundles as opposed to simply vector bundles. 
Part (c) follows from the choice of $F$.
\end{proof}

We use the proposition to define a grading on Higgs triples as follows.

\begin{defn}
For any stable Higgs triple $\tri$ such that $(\del, \theta)$ is unstable, we 
define the $\mu_+$ grading to be given by
$$ gr^+\tri = (gr^+(\del, \theta), \psi), $$
where $gr^+(\del, \theta)$ and $\psi$ are as above. If $(\del, \theta)$ is 
semistable,
$$ gr^+\tri = (Gr(\del, \theta), 0), $$
where $Gr(\del, \theta)$ is the usual HN filtration for Higgs bundles.
\end{defn}

The two gradings $gr^-$ and $gr^+$ will help us stratify the space 
$\hatb$. At this point the reader may notice that for any triple $\tri$, both 
gradings naturally yield $(n+2)$-tuples of integers $(d_\phi, r_\phi, r_1, 
\ldots, r_n)$ which satisfy all the conditions of Lemma \ref{tuples lemma}. The 
numbers are obtained from the gradings in the following manner: given a $gr^-$ 
grading, let $d_\phi$ and $r_\phi$ be the degree and rank of $E_\phi$ 
respectively, where $r_\phi$ is the rank of the bundle $F_i/F_{i-1}$. Given a 
$gr^+$ grading, let $d_\phi$ and $r_\phi$ be the degree and rank of the 
bundle $F_{n+1}/F_n$, where $r_\phi$ is the same as in $gr^-$. We now proceed to 
the stratification.

\begin{defn}
Given any $(n+2)$-tuple $(d_\phi, r_\phi, r_1, \ldots, r_n) \in 
\mathcal{I}_\tau$, let
\begin{multline*}
\mathcal{W}^{\pm}(d_\phi, r_\phi, r_1, \ldots, r_n) = \Big\{ \tri\in\hatb\, 
\Big| \, \tri \textrm{ is a stable 
Higgs triple, and } \\
gr^{\pm}\tri \in \mathcal{Z}(d_\phi, r_\phi, r_1, \ldots, r_n) \Big\} \bigcup 
\mathcal{Z}(d_\phi, r_\phi, r_1, \ldots, r_n).
\end{multline*}
and
$$ \mathcal{W}^{\pm}_\tau = \bigcup\limits_{(d_\phi, r_\phi, r_1, \ldots, r_n) 
\in \mathcal{I}_\tau} \mathcal{W}^{\pm}(d_\phi, r_\phi, r_1, \ldots, r_n) $$
\end{defn}

The subspaces $\mathcal{W}^+_\tau$ form a stratification of $\hatb$. A similar 
result is also true for $\mathcal{W}^-$. Before proving the next proposition, we 
introduce the following notation: let $\mathcal{U}o^{\pm}_\tau$ be the 
corresponding stratifying spaces of $\tau$-semistable pairs as in \cite{BDW}, 
and $\hatb o$ the master space of all such semistable pairs.

\begin{prop}
If $r > 2$, then for critical values $\tau \in \left(\frac{d}{r}, \frac{d}{r-1} 
\right)$, the complex codimension of $\mathcal{W}^{\pm}_\tau$ in $\hat{\B}$ is 
at least 1.
\end{prop}

\begin{proof}
We first consider $\mathcal{W}^+_\tau$. Let
$$ \mathcal{W}o^+_\tau = \left\{ (\del, \phi)\in \mathcal{A}\times H^0(E^{\del}) 
\,\Big| \, \tri \in \mathcal{W}^+_\tau \right\}. $$
Then $\mathcal{W}^+_\tau$ fibres over $\mathcal{W}o^+_\tau$ with fibre $\{\tri 
\in \mathcal{W}^+_\tau\,|\, \theta(\phi) = 0\}$ over the point $(\del, \phi) \in 
\mathcal{W}o^+_\tau$. Similarly, if we let
$$ \mathcal{U}^+_\tau = \left\{ \tri\in\mathcal{W}^+_\tau \, \Big| \, (\del, 
\phi) \textrm{ is } \tau\textrm{-semistable as a pair} \right\}, $$
then the space $\mathcal{U}^+_\tau$ fibres over $\mathcal{U}o^+_\tau$. However, 
the space $\mathcal{U}o^+_\tau$ of $\tau$-semi-stable pairs is 
an open subset of $\mathcal{W}o^+_\tau$. Moreover, the space 
$\mathcal{U}o^+_\tau$ has positive codimension in $\hatb o$, the master space of 
pairs. Since $\hatb$ fibres in a similar manner over $\hatb o$, we see that 
the space $\mathcal{U}^+_\tau$ must have positive codimension in $\hatb$, and 
hence so must $\mathcal{W}^+_\tau$. We summarize the above argument in the 
following diagram.
\begin{equation*}
\begin{CD}
\hatb @<<< \mathcal{U}^+_\tau @>>> \mathcal{W}^+_\tau \\
@VVV @VVV @VVV \\
\hatb o @<pos. codim.<< \mathcal{U}o^+_\tau @>open>> \mathcal{W}o^+_\tau \\
\end{CD}
\end{equation*}
The same argument works for the case $\mathcal{W}^-_\tau$.
\end{proof}

\section{Morse Theory and Birational Equivalence}
\setcounter{theorem}{0}
\setcounter{equation}{0}

We now look more carefully at the Morse theory of the function $f$. We will see 
that the stratification obtained using Morse theory coincides 
with the algebraic stratification of the previous section. We first define the 
flow on $\hatb$ as follows.

\begin{defn}
Let $\Phi: \hatb \times (-\infty, \infty) \rightarrow \hatb$ be the flow
$$ \Phi_t([\del, \theta, \phi]) = [\del, \theta, e^{-t/2\pi r}\phi]. $$
\end{defn}

Note that since the flow only affects the section $\phi$, our situation is 
identical to the one using stable pairs in \cite{BDW}.

\begin{prop}
$\Phi$ is continuous and it preserves $\hatb_0$. Further, it coincides with the 
gradient flow of $f$ on $\hatb_0$ i.e.,
$$ \frac{d\Phi_t}{dt} = -\nabla_{\Phi_t} f. $$
\end{prop}
\begin{proof}
This is similar to proposition 4.1 in \cite{BDW}. Since the tangent space 
$T_{[\del, \theta, \phi]}\hatb_0$ has been identified with $H^1({\bf 
C}_0^{\tri})$, the infinitesimal vector field of the $S^1$ action on $\hatb_0$ 
is given by $\xi^{\#}[\del, \theta, \phi] = \frac{i}{r}(0,0,r)$. Hence,
\begin{align*}
\nabla_{\Phi_t[\del,\theta,\phi]} f \quad &= \quad \frac{-1}{2\pi i} 
\nabla_{\Phi_t[\del,\theta,\phi]}\Psi \quad = \quad \frac{1}{2\pi i} 
\xi^{\#}(\Phi_t[\del, \theta, \phi]) \\
&= \quad \frac{1}{2\pi r}(0,0,e^{-t/2\pi r}\phi) \quad = \quad 
-\frac{d}{dt}(0,0,e^{-t/2\pi r}\phi) \\
&= \quad - \frac{d\Phi_t[\del,\theta,\phi]}{dt}.
\end{align*}
\end{proof}

\begin{defn}
Given a critical value $\tau$, and $(d_\phi, r_\phi, r_1, \ldots, r_n) \in 
\mathcal{I}_\tau$, let
$$\mathcal{W}^s(d_\phi, r_\phi, r_1, \ldots, r_n) = \left\{[\del, \theta, \phi] 
\in \hatb\, \Big| \, \lim\limits_{t\rightarrow\infty} \Phi_t([\del, \theta, 
\phi]) \in \mathcal{Z}(d_\phi, r_\phi, r_1, \ldots, r_n) \right\}. $$
We similarly define $\mathcal{W}^u(d_\phi, r_\phi, r_1, \ldots, r_n)$ by taking 
the limit as $t\rightarrow -\infty$. This allows us to define
$$ \mathcal{W}^s_\tau = \bigcup\limits_{(d_\phi, r_\phi, r_1, \ldots, r_n) \in 
\mathcal{I}_\tau} \mathcal{W}^s(d_\phi, r_\phi, r_1, \ldots, r_n) $$
and similarly,
$$ \mathcal{W}^u_\tau = \bigcup\limits_{(d_\phi, r_\phi, r_1, \ldots, r_n) \in 
\mathcal{I}_\tau} \mathcal{W}^u(d_\phi, r_\phi, r_1, \ldots, r_n). $$
We shall call $\mathcal{W}^s$ and $\mathcal{W}^u$ the stable and unstable Morse 
stratifications of $\hatb$, respectively.
\end{defn}

\begin{prop}
For every critical value $\tau$, the Morse stratification of $\hatb$ coincides 
exactly with the algebraic stratification of $\hatb$ i.e.,
$$ \mathcal{W}^s_\tau = \mathcal{W}^+_\tau \quad \textrm{ and } \quad 
\mathcal{W}^u_\tau = \mathcal{W}^-_\tau. $$
\end{prop}
\begin{proof}
This is similar to proposition 4.3 in \cite{BDW}. First consider 
$\mathcal{W}^u_\tau$ and $\mathcal{W}^-_\tau$. Since both stratify 
$\hatb$, it is sufficient to show that $\mathcal{W}^-_\tau \subseteq 
\mathcal{W}^u_\tau$ for all $\tau$. In particular, it is enough to show 
that
$$ \mathcal{W}^-_\tau(d_\phi, r_\phi, r_1, \ldots, r_n) \subseteq 
\mathcal{W}^u_\tau(d_\phi, r_\phi, r_1, \ldots, r_n) $$
for all $(d_\phi, r_\phi, r_1, \ldots, r_n) \in \mathcal{I}_\tau$. Consider 
any $[\del, \theta, \phi] \in \mathcal{W}^-_\tau(d_\phi, r_\phi, r_1, \ldots, 
r_n)$. Let
$$ 0 = E_\phi = F_0 \subset F_1 \subset \cdots \subset F_n = E $$
denote the $\mu_-$ filtration of $\tri$. Fix real numbers $0 < \mu_1 < \mu_2 < 
\cdots < \mu_n$ such that $\D\sum\limits_{i=1}^{n} r_i\mu_i = r_\phi$, and 
consider the 1-parameter subgroup of gauge tranformations in 
$\mathcal{G}_0^{\C}$,
$$ g_t = 
\begin{pmatrix}
e^{t/2\pi} & 0 & 0 & \cdots & 0 \\
0 & e^{-t\mu_1/2\pi r} & 0 & \cdots & 0 \\
\vdots & & \ddots & & \vdots \\
0 & \cdots & 0 & 0 & e^{-t\mu_n/2\pi r} \\
\end{pmatrix} $$
written with respect to the filtration above. Then
$$ \lim\limits_{t\rightarrow -\infty} \Phi_t[\del,\theta,\phi] \, = \, 
\lim\limits_{t\rightarrow -\infty} [\del,\theta,g_t^{-1}\phi] \, = \, 
\lim\limits_{t\rightarrow -\infty} [g_t(\del, \theta),\phi] = [gr^-(\del, 
\theta), \phi]. $$
The other case is similar. Again, it suffices to show that $$ 
\mathcal{W}^+_\tau(d_\phi, r_\phi, r_1, \ldots, r_n) \subseteq 
\mathcal{W}^s_\tau(d_\phi, r_\phi, r_1, \ldots, r_n) $$
for all $(d_\phi, r_\phi, r_1, \ldots, r_n) \in \mathcal{I}_\tau$. So take any 
$[\del, \theta, \phi] \in \mathcal{W}^+_\tau(d_\phi, r_\phi, r_1, \ldots, r_n)$. 
Let 
$$ 0 = F_0 \subset F_1 \subset \cdots \subset F_{n+1} = E $$
denote the $\mu_+$ filtration of $\tri$. If $(\del, \theta)$ is semistable (as a 
Higgs bundle), fix real numbers $1 > \mu_1 > \mu_2 > \cdots > \mu_{n+1}$ such 
that $\D\sum\limits_{i=1}^{n} r_i\mu_i = 0$, and let $g_t$ be the 
gauge tranformation 
$$ g_t = 
\begin{pmatrix}
e^{t\mu_1/2\pi r} & 0  & \cdots & 0 \\
\vdots & & \ddots & \vdots \\
0 & \cdots & 0 & e^{-t\mu_{n+1}/2\pi r} \\
\end{pmatrix} $$
written with respect to the filtration above. If $(\del, \theta)$ is not 
semistable (as a Higgs bundle), then let $\mu_{n+1} = 1$, and choose the other 
$\mu_i$ such that $\D r_{n+1} + \sum\limits_{i=1}^{n} r_i\mu_i = 0$. The rest of 
the argument proceeds as before.
\end{proof}

We finally focus our attention on how the spaces 
$\mathcal{B}_\tau$ are affected as $\tau$ is varied smoothly. We start with a 
proposition which follows essentially from the Morse theory of $f$.

\begin{theorem}
\begin{enumerate}
\item[(a)] If $f$ has no critical values in the interval $[\tau, 
\tau+\varepsilon]$, then the Morse flow induces a biholomorphism between 
$\B_{\tau+\varepsilon}$ and $\B_\tau$.
\item[(b)] If $\tau$ is the only critical value of $f$ in the interval $[\tau, 
\tau+\varepsilon]$, then the Morse flow induces a biholomorphism between 
$\B_{\tau+\varepsilon}\backslash \mathbb{P}_\varepsilon(\mathcal{W}^+_\tau)$ and 
$\B_\tau \backslash \mathcal{Z}_\tau$, where
$$ \mathbb{P}_\varepsilon(\mathcal{W}^+_\tau) = (\mathcal{W}^+_\tau \cap 
f^{-1}(\tau+\varepsilon))/U(1). $$
\end{enumerate}
\end{theorem}

\begin{proof}
This proof follows the same argument as the corresponding statement, 
theorem 4.4 in \cite{BDW}. As a result, we simply sketch here a construction of 
the biholomorphism.

Denote the equivalence class of the triple $[\del, \theta, \phi] \in \hatb$ by 
$x$. Let
$$ F: f^{-1}(\tau + \varepsilon) \times [0, \infty) \longrightarrow \R $$ 
be defined by $F(x, t) = f(\Phi_t(x))$. In part (a), since $f$ is by assumption, 
smooth on the 
interval $[\tau, \tau + \varepsilon]$ , we see that $F$ is smooth. 
Further,
$$ \frac{\partial F}{\partial t}\Big|_{(x, t)} = df_{\Phi_t(x)}\left( 
\frac{\partial \Phi_t}{\partial x}\Big|_x \right) = \parallel \nabla_{\Phi_t(x)} 
f \parallel^2 \neq 0, $$
for any $(x, t) \in F^{-1}(\tau)$ since $\tau = f(\Phi_t(x))$ is not a critical 
value of $f$. As a result, we may use the Implicit Function Theorem to solve the 
equation $F(x, t) = \tau$ and get $t = t(x)$ as a smooth function of $x$. This 
allows us to define a map:
\begin{equation}
\hat{\sigma}_+: f^{-1}(\tau + \varepsilon) \longrightarrow f^{-1}(\tau)
\end{equation}
by $\hat{\sigma}_+(x) = \Phi_{t(x)}(x)$.
$\hat{\sigma}_+$ then induces a biholomorphism
$$\sigma: \B_{\tau + \varepsilon} = f^{-1}(\tau + \varepsilon)/U(1) 
\longrightarrow f^{-1}(\tau)/U(1) = \B_\tau.$$
The above follows once we establish that the map $\hat{\sigma}_+$ and the 
complex structure on the spaces $f^{-1}(\tau + \varepsilon)$, $f^{-1}(\tau)$ are 
all $U(1)$-invariant. In part (b), the same argument as in part (a) gives a 
smooth map
$$ \hat{\sigma}_+: f^{-1}(\tau + \varepsilon)\backslash\mathcal{W}^+_{\tau}  
\longrightarrow f^{-1}(\tau)\backslash\mathcal{Z}_\tau. $$
This map can then be extended (continuously) across $\mathcal{W}^+_\tau$ by 
setting $\hat{\sigma}_+(x) = \lim\limits_{t\rightarrow \infty} \Phi_t(x)$ for $x 
\in \mathcal{W}^+_\tau$. Finally, the same argument as before shows that 
$\hat{\sigma}_+$ is a biholomorphism onto its image away from 
$\mathbb{P}_\varepsilon(\mathcal{W}^+_\tau) = (\mathcal{W}^+_\tau \cap 
f^{-1}(\tau+\varepsilon))/U(1)$.
\end{proof}

Note that if the flow lines are reversed, the above argument shows that a 
similar statement must be true for $\B_{\tau-\varepsilon}$ and $\B_\tau$. This 
observation leads us to the following.

\begin{picture}(350,60)(30,0)
\put(163,12){\line(1,-1){25}}
\put(188,-13){\line(-1,0){3}}
\put(188,-13){\line(0,1){3}}
\put(232,12){\line(-1,-1){25}}
\put(207,-13){\line(1,0){3}}
\put(207,-13){\line(0,1){3}}
\end{picture}
\vspace{-0.9in}

\begin{cor}
If $\tau$ is the only critical value in $[\tau - \varepsilon, \tau + 
\varepsilon]$, then $\B_{\tau-\varepsilon}$ and $\B_{\tau + \varepsilon}$ are 
related as shown:
$$
\begin{CD}
\B_{\tau - \varepsilon} \qquad\qquad\qquad 
\B_{\tau+\varepsilon} \\
\sigma_- \quad\qquad\qquad \sigma_+ \\
\B_\tau
\end{CD}
$$
where $\sigma_{\pm}$ are continuous, and
$\sigma_{\pm}: \B_{\tau \pm \varepsilon}\backslash \sigma_{\pm}^{-1}( 
\mathcal{Z}_\tau) \longrightarrow \B_\tau \backslash \mathcal{Z}_\tau$
are biholomorphisms.
\end{cor}

\begin{theorem}
For all noncritical values of $\tau$ in $\left(\frac{d}{r}, \frac{d}{r-1} 
\right)$, the spaces $\B_\tau$ are all birational.
\end{theorem}

\begin{proof}
By the previous corollary, the complex manifolds $\B_{\tau \pm 
\varepsilon}\backslash \sigma_{\pm}^{-1}(\mathcal{Z}_\tau)$ are biholomorphic. 
Let $\tau \in \left(\frac{d}{r}, \frac{d}{r-1} \right)$ be a critical value. 
Let
$$ T_- = \left\{(\del, \phi)\,\Big|\, \tri\in \B_{\tau - \varepsilon}\backslash 
\sigma_{-}^{-1}(\mathcal{Z}_\tau)\textrm{ for some Higgs field }\theta \right\} 
$$
$$ T_+ = \left\{(\del, \phi)\,\Big|\, \tri\in \B_{\tau + \varepsilon}\backslash 
\sigma_{+}^{-1}(\mathcal{Z}_\tau)\textrm{ for some Higgs field }\theta \right\} 
$$
Then $\B_{\tau - \varepsilon}\backslash \sigma_{-}^{-1}(\mathcal{Z}_\tau)$ is a 
fibration over $T_-$ and $\B_{\tau + \varepsilon}\backslash 
\sigma_{+}^{-1}(\mathcal{Z}_\tau)$ is a fibration over $T_+$. However, the sets 
of $(\tau \pm \varepsilon)$-stable Bradlow pairs are open subsets of $T_{\pm}$. 
Since these spaces of Bradlow pairs are birational from \cite{BDW}, 
it follows that $T_-$ is birational to $T_+$.  Finally, our biholomorphism is an 
extension by identity (as 
$\theta$ is untouched by the Morse flow) of the map $T_- \rightarrow T_+$ to 
the fibres of $\B_{\tau \pm \varepsilon}\backslash 
\sigma_{\pm}^{-1}(\mathcal{Z}_\tau)$ over $T_{\pm}$. It follows that
$\B_{\tau \pm \varepsilon}$ are birational. This completes the proof.
\end{proof}

Now we recall that $\B_{\frac{d}{r}}$ and $\B_{\frac{d}{r-1}}$ are moduli
spaces of Higgs bundles of rank $r$ and $r-1$, respectively. Let 
$$ \HH_0(r,d) = \left\{ (\del, \theta) \in \HH(r,d)\,
\Big|\, \textrm{det}(\theta) = 0 \right\} $$
We see that $\B_{\frac{d}{r} + \varepsilon}$ is a fibration
over $\HH_0$. We shall denote the resulting fibre space by $\mathcal{F}$. At the
other end we have a similar situation: we see that
$\B_{\frac{d}{r-1}-\varepsilon}$ is a projective bundle over
$\mathcal{M}_{Higgs}(r-1,d)$ with fibre over $(E', \theta)$, the
projectivization of Higgs extensions of $(E', \theta)$ by $ (\mathcal{O},$
const.$)$. Denote this bundle by $\mathcal{E}$. This proves our second result, 
and gives us a geometric relationship between the moduli spaces of Higgs bundles 
of rank $r$ and $r-1$.

\begin{theorem}
The fibre space $\mathcal{F}$ over $\HH_0(r,d)$ is birational to the projective 
bundle $\mathcal{E}$ over $\HH(r-1,d)$.
\end{theorem}

Although we have successfully established a geometric correspondence between 
moduli spaces of Higgs bundles, some questions still remain unanswered. For 
instance, an explicit description of the space $\mathcal{E}$ and $\mathcal{F}$ 
would be useful. In general, the space $\mathcal{F}$ appears to have rather 
badly 
behaved fibres. However, we conjecture that these fibres have constant dimension 
over the nilpotent cone, which sits inside $\HH_0$. 
Another motivation in obtaining this result is to understand the torus
action on the moduli space of Higgs bundles using induction on rank. We hope to 
study this in the near future as well.

{\it Acknowledgements}.
The results in this paper were submitted as part of my Ph.D. dissertation in the department of Mathematics, University of Chicago. I would like to thank my advisor Kevin Corlette for his encouragement and guidance. I am also grateful to Steve Bradlow, Vladimir Baranovsky, and Madhav Nori for several useful discussions.

\end{document}